\documentclass{elsart3-1}
\usepackage{amssymb,amsmath,amsfonts,latexsym}
\usepackage[english,francais]{babel}

\newtheorem{theorem}{Theorem}[section]
\newtheorem{lemma}[theorem]{Lemma}
\newtheorem{proposition}[theorem]{Proposition}
\newtheorem{corollary}[theorem]{Corollary}
\newtheorem{definition}[theorem]{Definition\rm}
\newtheorem{remark}{\it Remark\/}
\newtheorem{example}{\it Example\/}
\renewcommand{\theequation}{\arabic{equation}}
\def\og{\leavevmode\raise.3ex\hbox{$\scriptscriptstyle\langle\!\langle$~}}
\def\fg{\leavevmode\raise.3ex\hbox{~$\!\scriptscriptstyle\,\rangle\!\rangle$}}

\newcommand{\nc}[2]{\newcommand{#1}{#2}}
\newcommand{\rnc}[2]{\renewcommand{#1}{#2}}
\rnc{\theequation}{\thesection.\arabic{equation}}
\def\note#1{{}}
\def\Label#1{\label{#1}\ifmmode\llap{[#1] }\else
\marginpar{\smash{\hbox{[#1]}}}\fi}

\nc{\beq}{\begin{equation}}
\nc{\eeq}{\end{equation}}
\rnc{\[}{\beq}
\rnc{\]}{\eeq}
\rnc{\Label}{\label}
\nc{\ba}{\begin{array}}
\nc{\ea}{\end{array}}
\nc{\bea}{\begin{eqnarray}}
\nc{\beas}{\begin{eqnarray*}}
\nc{\eeas}{\end{eqnarray*}}
\nc{\eea}{\end{eqnarray}}
\nc{\be}{\begin{enumerate}}
\nc{\ee}{\end{enumerate}}
\nc{\bd}{\begin{diagram}}
\nc{\ed}{\end{diagram}}
\nc{\bi}{\begin{itemize}}
\nc{\ei}{\end{itemize}}
\nc{\bpr}{\begin{proposition}}

\nc{\bth}{\begin{theorem}}
\nc{\ble}{\begin{lemma}}
\nc{\bco}{\begin{corollary}}
\nc{\bre}{\begin{remark}}
\nc{\bex}{\begin{example}}
\nc{\bde}{\begin{definition}}
\nc{\ede}{\end{definition}}
\nc{\epr}{\end{proposition}}
\nc{\ethe}{\end{theorem}}
\nc{\ele}{\end{lemma}}
\nc{\eco}{\end{corollary}}
\nc{\ere}{\hfill\mbox{$\losenge$}\end{remark}}
\nc{\eex}{\hfill\mbox{$\losenge$}\end{example}}
\nc{\bpf}{{\it Proof.~~}}
\nc{\epf}{\hfill\mbox{$\square$}\vspace*{3mm}}
\nc{\hsp}{\hspace*}
\nc{\vsp}{\vspace*}

\def\ot{\otimes}
\nc{\te}{\!\ot\!}
\nc{\bmlp}{\mbox{\boldmath$\left(\right.$}}
\nc{\bmrp}{\mbox{\boldmath$\left.\right)$}}
\nc{\LAblp}{\mbox{\LARGE\boldmath$($}}
\nc{\LAbrp}{\mbox{\LARGE\boldmath$)$}}
\nc{\Lblp}{\mbox{\Large\boldmath$($}}
\nc{\Lbrp}{\mbox{\Large\boldmath$)$}}
\nc{\lblp}{\mbox{\large\boldmath$($}}
\nc{\lbrp}{\mbox{\large\boldmath$)$}}
\nc{\blp}{\mbox{\boldmath$($}}
\nc{\brp}{\mbox{\boldmath$)$}}
\nc{\LAlp}{\mbox{\LARGE $($}}
\nc{\LArp}{\mbox{\LARGE $)$}}
\nc{\Llp}{\mbox{\Large $($}}
\nc{\Lrp}{\mbox{\Large $)$}}
\nc{\llp}{\mbox{\large $($}}
\nc{\lrp}{\mbox{\large $)$}}
\nc{\lbc}{\mbox{\Large\boldmath$,$}}
\nc{\lc}{\mbox{\Large$,$}}
\nc{\Lall}{\mbox{\Large$\forall\;$}}
\nc{\bc}{\mbox{\boldmath$,$}}
\nc{\ra}{\rightarrow}
\nc{\ci}{\circ}
\nc{\cc}{\!\ci\!}
\nc{\lra}{\longrightarrow}
\nc{\imp}{\Rightarrow}
\rnc{\iff}{\Leftrightarrow}
\nc{\inc}{\mbox{$\,\subseteq\;$}}
\rnc{\subset}{\inc}
\nc{\0}{\sb{(0)}}
\nc{\1}{\sb{(1)}}
\nc{\2}{\sb{(2)}}
\nc{\3}{\sp{(3)}}
\nc{\4}{\sp{(4)}}
\nc{\5}{\sp{(5)}}
\nc{\6}{\sp{(6)}}
\nc{\7}{\sp{(7)}}

\def\tr{{\rm tr}}
\def\Tr{{\rm Tr}}
\def\st{\stackrel}
\def\<{\langle}
\def\>{\rangle}

\def\ker{\mbox{$\mathop{\mbox{\rm Ker$\,$}}$}}

\def\aut{\mbox{$\mathop{\mbox{\rm Aut}}$\,}}

\nc{\spp}{\mbox{${\cal S}{\cal P}(P)$}}
\nc{\ob}{\mbox{$\Omega\sp{1}\! (\! B)$}}
\nc{\op}{\mbox{$\Omega\sp{1}\! (\! P)$}}
\nc{\oa}{\mbox{$\Omega\sp{1}\! (\! A)$}}
\nc{\dr}{\mbox{$\Delta_{R}$}}
\nc{\dsr}{\mbox{$\Delta_{\Omega^1P}$}}
\nc{\ad}{\mbox{$\mathop{\mbox{\rm Ad}}_R$}}
\nc{\as}{\mbox{$A(S^3\sb s)$}}
\nc{\bs}{\mbox{$A(S^2\sb s)$}}
\nc{\slc}{\mbox{$A(SL(2,\C))$}}
\nc{\suq}{\mbox{$\cO(SU_q(2))$}}
\nc{\tc}{\widetilde{can}}

\rnc{\epsilon}{\varepsilon}
\rnc{\phi}{\varphi}
\nc{\ha}{\mbox{$\alpha$}}
\nc{\hb}{\mbox{$\beta$}}
\nc{\hg}{\mbox{$\gamma$}}
\nc{\hd}{\mbox{$\delta$}}
\nc{\he}{\mbox{$\varepsilon$}}
\nc{\hz}{\mbox{$\zeta$}}
\nc{\hs}{\mbox{$\sigma$}}
\nc{\hk}{\mbox{$\kappa$}}
\nc{\hm}{\mbox{$\mu$}}
\nc{\hn}{\mbox{$\nu$}}
\nc{\hl}{\mbox{$\lambda$}}
\nc{\hG}{\mbox{$\Gamma$}}
\nc{\hD}{\mbox{$\Delta$}}
\nc{\hT}{\mbox{$\Theta$}}
\nc{\ho}{\mbox{$\omega$}}
\nc{\hO}{\mbox{$\Omega$}}
\nc{\hp}{\mbox{$\pi$}}
\nc{\hP}{\mbox{$\Pi$}}

\nc{\qpb}{quantum principal bundle}

\def\C{{\Bbb C}}
\def\N{{\Bbb N}}

\def\Z{{\Bbb Z}}

\def\cO{{\mathcal O}}

\def\cT{{\mathcal T}}
\def\cK{{\mathcal K}}

\def\nc{\newcommand}

\def\tc{\tilde c}

\def\st{\stackrel}

\newcommand{\fa}{\forall}

\def\ra{\rightarrow}

\def\te{{\tilde e}}

\begin{document}
\begin{frontmatter}
\selectlanguage{english}
\vspace*{-95pt}
\title{Noncommutative index theory for mirror quantum spheres}

\vspace{-2.6cm}

\selectlanguage{francais}
\title{}

\selectlanguage{english}
\author[authorlabel1,authorlabel2]{Piotr M.\ Hajac},
\ead{http://www.fuw.edu.pl/$\!\widetilde{\phantom{m}}\!$pmh}
\author[authorlabel2]{Rainer Matthes},
\ead{matthes@fuw.edu.pl}
\author[authorlabel3]{Wojciech Szyma\'nski}
\ead{Wojciech.Szymanski@newcastle.edu.au}
\address[authorlabel1]{Instytut Matematyczny, Polska Akademia Nauk,
ul.\ \'Sniadeckich 8, Warszawa, 00--956 Poland}
\address[authorlabel2]{Katedra Metod Matematycznych Fizyki, Uniwersytet
Warszawski, ul.\ Ho\.za 74, Warszawa, 00--682 Poland}
\address[authorlabel3]{School of Mathematical and Physical Sciences,
University of Newcastle, Callaghan, NSW 2308, Australia}

\begin{abstract}
\vskip 0.5\baselineskip

\noindent We introduce and analyse a new type of quantum 2-spheres.
Then we apply  index theory for noncommutative line bundles
over these spheres to conclude that quantum lens spaces are 
non-crossed-product examples
of principal extensions of $C^*$-algebras.
{\it To cite this article:  P.M.~Hajac, R.~Matthes, W.~Szyma\'nski}

\selectlanguage{francais}
\vskip 0.5\baselineskip
\noindent{\bf R\'esum\'e}
\vskip 0.5\baselineskip
\noindent{\bf Th\'eorie de l'indice non commutative pour des sph\`eres quantiques miroirs.}
Nous introduisons et analysons un nouveau type de 2-sph\`eres quantiques. 
Nous appliquons la th\'eorie de l'indice pour les fibr\'es lin\'eaires non 
commutatifs sur ces sph\`eres afin de d\'eduire que les espaces de lentilles 
quantiques sont des exemples d'extensions principales de $C^*$-alg\`ebres
qui ne sont pas des produits crois\'es. 
%
{\it Pour citer cet article: P.M.~Hajac, R.~Matthes, W.~Szyma\'nski}
\end{abstract}
\end{frontmatter}

\noindent{\bf Version fran\c caise abr\'eg\'ee}
\vspace*{3mm}

\noindent
La $C^*$-alg\`ebre $C(S^2_{pq+})$ d'une sph\`ere quantique de Podle\'s
g\'en\'erique \cite{p-p87} peut \^etre vu comme un produit fibr\'e de deux alg\`ebres de Toeplitz $\cT$ \cite{s-ajl91}.
Le produit fibr\'e est d\'efini par la fonction symbole
$\cT\st{\sigma}{\ra}C(S^1)$ de la fa\c con suivante :
$C(S^2_{pq+})\cong\{(x,y)\in\cT\times\cT\,|\,\hs(x)=\hs(y)\}$. 
Cela peut \^etre vu comme un collage de deux disques quantiques dont on identifie les bords $S^1$.
%
Nous consid\'erons un collage similaire de disques quantiques mais avec une identification diff\'erente des bords, plus pr\'ecisement
$C(S^2_{pq-}):=\{(x,y)\in\cT\times\cT\,|\,\hs(x)=\overline{\hs}(y)\}$
(o\`u $\overline{\hs}$ est le morphisme d'alg\`ebre d\'efini en envoyant l'isom\'etrie g\'en\'eratrice de $\cT$ sur {\em l'inverse} de 
l'unitaire g\'en\'eratrice de $C(S^1)$). 
%
Cela peut \^etre vu comme le collage du bord d'un disque quantique avec celui son miroir image.
Du point de vue alg\'ebrique, nous commen\c cons par $\cO(D_p)$, la 
sous-$*$-alg\`ebre dense de $\cT$ d\'efinie par la $*$-alg\`ebre unitaire 
universelle pour la relation $x^*x-pxx^*=1-p$, $p\in [0,1)$, \cite{kl93}.
%
Nous consid\'erons alors les versions alg\'ebriques de 
$\hs$ et de $\overline{\hs}$ qui envoient $x$ respectivement
vers le g\'en\'erateur unitaire des polyn\^omes de Laurent et
vers son inverse. Nous \'etudions le produit fibr\'e 
$\{(x,y)\in\cO(D_p)\times\cO(D_q)\,|\,\hs(x)=\overline{\hs}(y)\}$ 
qui se trouve \^etre \'equivalent \`a la $*$-alg\`ebre unitaire universelle
g\'en\'er\'ee par $C$, $E$ et $F$ munis des relations
$
C^*C=1-pE-F,\quad CC^*=1-E-qF,\quad
EC=pCE,\quad CF=qFC,\quad
EF=0,\quad E=E^*,\quad F=F^*,\quad p,q\in [0,1]
$.
Nous l'appelons l'alg\`ebre des coordon\'ees d'une {\em 2-sph\`ere quantique miroir} et nous la notons par $\cO(S^2_{pq-})$.
Puisque l'alg\`ebre de Toeplitz $\cT$ est la $C^*$-alg\'ebre enveloppante de 
$\cO(D_p)$ \cite{kl93}, nous avons que 
$C(S^2_{pq-})$ est la $C^*$-alg\`ebre enveloppante de 
$\cO(S^2_{pq-})$.
%

La $C^*$-alg\`ebre $C(S^2_{pq-})$ a, \`a equivalence unitaire pr\`es, deux  
$*$-repr\'esentations born\'ees, irr\'eductibles, infini-dimensionnelles 
$\rho_+$ et $\rho_-$ et aussi une famille $\rho_\lambda$, ${\lambda\in S^1}$
de repr\'esentations uni-dimensionnelles. 
%
\note{Les noyaux de $\rho_+$ et de $\rho_-$ sont isomophiques aux op\'erateurs compacts $\cK$ et ont une intersection triviale. Cela donne la suite exacte
$
0 \rightarrow \cK\oplus\cK \rightarrow C(S^2_{pq-})
\rightarrow C(S^1) \rightarrow 0
$.
De l\`a, nous d\'eduisons les $K$-groupes :
$K_0(C(S^2_{pq-}))\cong\Z\oplus\Z$ et $K_1(C(S^2_{pq-}))=0$.
Les repr\'esentations donnent aussi l'espace des id\'eaux primitifs de
$C(S^2_{pq-})$ qui consistent en une paire de points dense et un cercle. Ceci est en accord avec les propri\'et\'es respectives de $C(S^2_{pq+})$.
}
Cette classification des repr\'esentations et beaucoup de propri\'et\'es de 
$C(S^2_{pq-})$ sont en accord avec celles de Podle\'s. 
N\'eanmoins, en consid\'erant les comportements diff\'erents des  
$K_0$-classes de deux projections minimales non-\'equivalentes dans
$C(S^2_{pq-})$ et  dans $C(S^2_{pq+})$, cela permet de montrer que les 
$C^*$-alg\`ebres $C(S^2_{pq+})$ et $C(S^2_{pq+})$ en sont pas Morita-\'equivalentes.
%
La diff\'erence entre ces deux alg\`ebres est aussi visible dans leur structure de module de Fredholm. Contrairement \`a $C(S^2_{pq+})$, la paire de repr\'esentations irr\'eductible infini-dimensionnelles 
$(\rho_+,\rho_-)$ n'est pas un module de Fredholm. A la place, nous obtenons :
%
\vspace*{1mm}
\ble \label{fFred-}
Soit $\Tr$ la trace. La paire 
$(\,\rho_+\oplus\rho_-\,,\,\int_{S^1}\rho_{\lambda}d\lambda\,)$ est un
module de Fredholm 1-sommable sur $\cO(S^2_{pq-})$, i.e.\
$
\mbox{$\fa$}\; x\in \cO(S^2_{pq-}):\;
\Tr\;\mbox{$\left|\,\rho_+(x)\oplus\rho_-(x)-\int_{S^1}\rho_{\lambda}(x)
d\lambda\,\right|$}\;
<\infty
$.
\ele
\vspace{3mm}
%

Soit $\cO(S^3_{pq\theta})$ l'alg\`ebre des coordonn\'ees de la 3-sph\`ere quantique de type  Heegaard \cite{bhms06}. On peut montrer que 
$\cO(S^2_{pq-})\subset \cO(S^3_{pq\theta})$ est une extension principale, dans le sens de  \cite{bh04}, pour l'action naturelle $U(1)$ sur $\cO(S^3_{pq\theta})$ (de cette fa\c con, nous pouvons voir cette extension comme une d\'eformation non commutative de la fibration de Hopf). Les  
$\cO(S^2_{pq-})$-modules associ\'es \`a cette extension, via des repr\'esentations uni-dimensionelles (fibr\'es lin\'eaires non commutatifs) peuvent \^etre d\'efinis comme  $\cO(S^3_{pq\theta})_{\mu}:=\{x\in\cO(S^3_{pq\theta})~|~\ha_g(x)=g^{-\mu}x,
\forall g\in U(1)\},~\mu\in\Z$. 
En appliquant \cite[Theorem 3.1]{bh04}, on peut montrer que les modules 
$\cO(S^3_{pq\theta})_{\mu}$ sont projectifs, finiment engendr\'es, et
en outre on peut 
d\'eterminer explicitement les idempotents $E_\mu$ qui les repr\'esentent.
En rempla\c cant $\cO(S^3_{pq\theta})$ par la $C^*$-alg\`ebre enveloppante
$C(S^3_{pq\theta})$ et en \'etendant par continuit\'e l'action $U(1)$,
nous obtenons une extension principale  
$C(S^2_{pq-})\subset C(S^3_{pq\theta})$ dans le sens de \cite{e-da00}, et
la famille $C(S^3_{pq\theta})_{\mu}$, $\mu\in\Z$ de $C(S^2_{pq-})$-modules.
On peut montrer que ces modules peuvent \^etre repr\'esent\'es par 
les m\^emes idempotents $E_\mu$.
Maintenant, nous pouvons repr\'esenter le module de Fredholm du Lemme~\ref{fFred-} par une trace $\tr$ sur $\cO(S^2_{pq-})$ et, en se servant de
\cite{h-pm00}, nous calculons les couplages de ces idempotents avec le module de Fredholm qui sont  les indices d'op\'erateurs de Fredholm.
\vspace*{3mm}\bth\label{findex}
Soit $\langle~,~\rangle$ le couplage de cohomologie cyclique et  
$K$-th\'eorie. Alors
$
\langle \tr,[E_\mu]\rangle = \mu$, $\fa\mu\in\Z
$.
\ethe
\vspace*{3mm}

Nous appliquons ce r\'esultat aux espaces de lentilles quantiques. Puisque 
$Z_n:=\Z/n\Z$ est un sous-groupe de $U(1)$, nous devons d\'efinir 
de nouvelles sous-alg\`ebres constitu\'ees de points fixes
$C(L^n_{pq\theta}):=C(S^3_{pq\theta})^{Z_n}$.
Elles jouent le r\^ole de l'alg\`ebre des fonctions continues sur les 
lentilles quantiques (cf.\ \cite{mt92,hs03}).
Comme dans le cas classique, il y a toujours une action naturelle de 
$U(1)$ sur $C(L^n_{pq\theta})$. On peut montrer que, pour tout 
$n\in\N\setminus\{0\}$, cette action est principale (libre) dans 
le sens de \cite{e-da00}, 
et que les sous-alg\`ebres constitu\'ees de points fixes 
$C(L^n_{pq\theta})^{U(1)}$ co\"\i ncident avec  $C(S^2_{pq-})$. 
En utilisant ces nouvelles actions $U(1)$, nous d\'efinissons 
les $C(S^2_{pq-})$-modules $C(L^n_{pq\theta})_\mu$, et on montre que
$C(L^n_{pq\theta})_\mu\cong C(S^3_{pq\theta})_{n\mu}\,$, $\fa\,n\in\N$.
D'un autre c\^ot\'e, nous montrons que les modules associ\'es avec 
les produits crois\'es par $\Z$ sont toujours libres.
En combinant ces r\'esultats avec le Th\'eor\`eme~\ref{findex}, 
nous obtenons notre r\'esultat principal :
\vspace*{3mm}\bth
Soient $n\in \N$, $n>0$, $p,q,\theta\in [0,1)$. 
La $C^*$-alg\`ebre $C(L^n_{pq\theta})$ est une $U(1)$-extension principale
de $C(S^2_{pq-})$, mais
$C(L^n_{pq\theta})\not\cong C(S^2_{pq-})\!\rtimes\!\Z$ comme $C^*$-alg\`ebres
pour toutes actions de $\Z$.
\ethe

\selectlanguage{english}

\section{Introduction}

\noindent
The purpose of this paper is twofold. First,
we present and study quantum two-spheres of  a new  topological type  (cf.~\cite{d-l03}). 
They are like
 twin siblings of the generic Podle\'s quantum spheres \cite{p-p87}, but of opposite
gender.
These new spheres were encountered while trying to understand why
the locally trivial quantum Hopf fibration built over a generic Podle\'s
sphere comes from an  action of $U(1)$ on $S^3$ that is equivalent to the
standard Hopf-fibration action by an orientation-reversing homeomorphism
\cite{hms06}. Taking the natural action leads 
 to a new fixed-point subalgebra
that coincides with the fibre product of Toeplitz algebras over $C(S^1)$
given by the symbol map and its composition with the automorphism
 inverting the unitary generator of $C(S^1)$. In contrast, the $C^*$-algebra of a 
 generic Podle\'s sphere 
is a fibre product of two Toeplitz 
algebras via two symbol maps (no inversion) \cite{s-ajl91}. Geometrically,
this means that the new sphere arises by gluing over the boundary
a quantum disc with its mirror image. It is precisely this
type of gluing of {\em oriented}  discs that yields an 
{\em oriented}  sphere. While in this classical case the orientation issue
does not show up in the $C^*$-algebra of continuous functions on $S^2$,
in the noncommutative setting it leads to two $C^*$-algebras that are
not even Morita equivalent. Since orientation is vital for the Dirac operator,
our new example may be a useful testing ground
for spectral triples (cf.\ \cite{dlps05}). Another curious aspect of
the new spheres is that, unlike the aforementioned Podle\'s
$C^*$-algebra, their $C^*$-algebra is not a graph algebra. However, it is very
close to  graph algebras \cite{hs02} and points at a more general concept of a graph algebra.  

\vspace*{3mm}

Our other main goal is to construct non-crossed-product examples of
principal extensions of $C^*$-algebras \cite{e-da00}. Much as there is more
to principal bundles than Cartesian products of groups with spaces, we
argue that there is  more to principal extensions of $C^*$-algebras than
crossed products. In the algebraic context of Hopf--Galois extensions, one
might prove that a given extension is not a crossed product
by arguing that algebraic crossed products have many
nontrivial invertible elements, whereas a given polynomial algebra
might lack such elements (e.g., see \cite[Appendix]{hm99}). This
argument does not apply on the $C^*$-level because it is easy to construct
nontrivial invertible elements in a $C^*$-algebra. It is difficult to
decide whether a given $C^*$-algebra is a crossed product with its 
fixed-point 
subalgebra because, contrary to the more restricted setting of von Neumann 
algebras, there appears to be no theory to manage such problems. Von Neumann 
algebras are at an opposite extreme to polynomial
algebras. For them theory tells us when an extension of algebras
has a crossed-product structure. (Think of measurable global sections of nontrivial
topological bundles.) Therefore, we find the $C^*$-context most interesting
to analyse this crossed-product problem. Our method follows
the classical pattern of algebraic topology: to prove that a principal
bundle is nontrivial, one computes an invariant of an associated vector
bundle. Herein we use the result of an index computation for finitely
generated projective modules over the $C^*$-algebra of a mirror quantum
sphere to prove that the $C^*$-algebras of noncommutative lens spaces
cannot be $\Z$-crossed products with the sphere  $C^*$-algebra. This
makes them essentially
 different from the standard example of the noncommutative
torus: 
$C(S^1)\cong C(T^2_\theta)^{U(1)}\inc C(T^2_\theta)\cong 
C(S^1)\rtimes_\theta\Z$.

\section{The comparison of Podle\'s and mirror quantum spheres}

\bde
Let $p,q\in[0,1]$. We call the universal unital $*$-algebra generated by
$C$, $E$ and $F$ satisfying
$
C^*C=1-pE-F,\quad CC^*=1-E-qF,\quad
EC=pCE,\quad CF=qFC,\quad
EF=0,\quad E=E^*,\quad F=F^*,\label{DE}
$
the coordinate algebra of a {\em mirror quantum 2-sphere}, and
denote it by $\cO(S^2_{pq-})$.
\ede
\vspace*{3mm}\noindent
The first two relations play the role of the sphere
equation. Here $E$ should be thought of as the coordinate measuring
the square of the distance of a point on the upper hemisphere from the equator plane,
while $F$ measures the square of the same distance but of a point on
the lower hemisphere. Indeed, if $p=1=q$, then we have the classical sphere,
$E$ and $F$ are functions with disjoint support, and the relationship
with the Cartesian coordinates is
$C=x+iy$, $\sqrt{E}-\sqrt{F}=z$. 
These new
spheres were found when examining a natural $U(1)$-action
on  Heegaard-type quantum 3-spheres $S^3_{pq\theta}$ \cite{bhms06}.
Recall that the polynomial algebra
 $\cO(S^3_{pq\theta})$ is the universal unital $*$-algebra
generated by elements $a$ and $b$  satisfying the  relations:
$
(1-aa^*)(1-bb^*)=0,\quad
a^*a=paa^*+1-p,\quad b^*b=qbb^*+1-q,\quad
ab=e^{2\pi i\theta} ba, \quad ab^*=e^{-2\pi i\theta}b^*a.
$
Here $p,q$ and $\theta$ are parameters with values in $[0,1)$.
Rescaling the generators by a unitary complex number,
i.e.\ $a\mapsto e^{i\phi}a$, $b\mapsto e^{i\phi}b$, gives an action
$\ha: U(1)\ra \aut\cO(S^3_{pq\theta})$. Moreover, one can show that 
the subalgebra of all \ha-invariants $\cO(S^3_{pq\theta})^{U(1)}$ is
independent of $\theta$ and isomorphic with
the unital $*$-algebra  $\cO(S^2_{pq-})$.

\vspace*{3mm}
The above defined $U(1)$-action extends by continuity to the
$C^*$-completion  $C(S^3_{pq\theta})$ of $\cO(S^3_{pq\theta})$. 
The fixed-point $C^*$-algebra
$C(S^3_{pq\theta})^{U(1)}$ turns out to be the universal (enveloping)
$C^*$-algebra
$C(S^2_{pq-})$ of the coordinate algebra of  $S^2_{pq-}$. Since
$C(S^3_{pq\theta})$ and  $C(S^3_{00\theta})$ are isomorphic via
a $U(1)$-equivariant map \cite[Theorem~2.8]{bhms06}, we can conclude
that $C(S^2_{pq-})\cong C(S^2_{00-})$ as $C^*$-algebras.
Another
way to think of this $C^*$-algebra is in terms of the fibre product of
two Toeplitz algebras given by the symbol map $\cT\st{\sigma}{\ra}C(S^1)$
and its composition with the automorphism of $C(S^1)$ sending the generating
unitary to its inverse:
$
C(S^2_{pq-})\cong\{(x,y)\in\cT\times\cT\,|\,\hs(x)=\overline{\hs}(y)\}
$.

\vspace*{3mm}\bpr \label{irrrep}
Let $\{\xi_k\}_{k\in\N}$ be an orthonormal basis of a Hilbert space,
and $\hl\in U(1)$.
Any irreducible bounded $*$-representation of
$\cO(S^2_{pq-})$
is unitarily equivalent to one of the following representations:
$
\rho_+(C)\xi_k=\sqrt{1-p^{k+1}}\,\xi_{k+1},\quad \rho_+(E)\xi_k=p^k\xi_k,
\quad \rho_+(F)\xi_k=0;\quad
\rho_-(C)\xi_k=\sqrt{1-q^k}\,\xi_{k-1},\quad
\linebreak
 \rho_-(E)\xi_k=0,
\quad \rho_-(F)\xi_k=q^k\xi_k$; or\quad
$\rho_{\lambda}(C)=\hl,
\quad \rho_{\lambda}(E)=0,
\quad \rho_{\lambda}(F)=0
$.
\note{
\begin{align}
&\rho_+(C)\xi_k=\sqrt{1-p^{k+1}}\,\xi_{k+1},\quad &&\rho_+(E)\xi_k=p^k\xi_k,
\quad &&\rho_+(F)\xi_k=0,
\label{rho+}\\
&\rho_-(C)\xi_k=\sqrt{1-q^k}\,\xi_{k-1},\quad &&\rho_-(E)\xi_k=0,
\quad &&\rho_-(F)\xi_k=q^k\xi_k,
\label{rho-}
\\ \label{1dim}
&\rho_{\lambda}(C)=\hl,
\quad &&\rho_{\lambda}(E)=0,
\quad &&\rho_{\lambda}(F)=0.
\end{align}
}
\epr
\vspace*{3mm}\noindent
We use the same symbols $\rho_+$, $\rho_-$ and $\rho_{\lambda}$ to denote
the extensions of these representations to the enveloping $C^*$-algebra
$C(S^2_{pq-})$. It turns out that the kernels of $\rho_+$ and $\rho_-$
have zero intersection and are both isomorphic to the algebra $\cK$ of compact operators
on a separable Hilbert space. 
The quotient of
$C(S^2_{pq-})$ by $\ker \rho_+ + \ker \rho_-$ is isomorphic to the algebra
of continuous
functions on the circle. Thus, the $C^*$-algebra $C(S^2_{pq-})$ admits
an exact sequence
$
0 \rightarrow \cK\oplus\cK \rightarrow C(S^2_{pq-})
\rightarrow C(S^1) \rightarrow 0.
$
Using this sequence or the fibre-product formula, one can compute:
$K_0(C(S^2_{pq-}))\cong\Z\oplus\Z$ and $K_1(C(S^2_{pq-}))=0$.

\vspace*{3mm}
The  properties of $C(S^2_{pq-})$ described so far match the respective properties
of the $C^*$-algebra of a generic Podle\'s sphere. One
can show that the latter is isomorphic with $C(S^2_{pq+})$, 
that is the universal unital
$C^*$-algebra for the relations
$\;
D^*D=1-pP-qQ,\quad DD^*=1-P-Q,
\quad
PD=pDP,\quad QD=qDQ,
\quad
PQ=0,\quad P=P^*,\quad Q=Q^*
$,
where $p,q\in [0,1)$. These relations appear almost identical to the defining
relations for $C(S^2_{pq-})$. Furthermore,  for the action given by
$a\mapsto e^{i\phi}a$, $b\mapsto e^{-i\phi}b$, we have
$C(S^2_{pq+})\cong C(S^3_{pq\theta})^{U(1)}$.
The fibre-product formula for $C(S^2_{pq+})$ also
differs only slightly from its counterpart for $C(S^2_{pq-})$, namely
$
C(S^2_{pq+})\cong\{(x,y)\in\cT\times\cT\,|\,\hs(x)=\hs(y)\}
$.
On the other hand, consider two inequivalent minimal projections
$e^\pm_1$ and $e^\pm_2$  in
$C(S^2_{pq\pm})$: it can be shown that their classes
are non-zero in the respective $K_0$-groups. However,
$[e^+_1]=-[e^+_2]$ in $K_0(C(S^2_{pq+}))$, while $[e^-_1]=[e^-_2]$
in $K_0(C(S^2_{pq-}))$. This leads to:
\vspace*{3mm}
\bpr \label{notiso}
The $C^*$-algebras $C(S^2_{pq+})$ and $C(S^2_{pq-})$ are not
stably isomorphic.
\epr
\vspace*{3mm}\noindent
Finally, recall that $C(S^2_{pq+})$ is a graph $C^*$-algebra
\cite{hs02}, whose primitive ideal space consists of a pair of points (which 
is dense in the entire space) and a circle. The primitive ideal space of 
$C(S^2_{pq-})$ turns out to be homeomorphic to that of $C(S^2_{pq+})$, 
but a careful analysis of all graph algebras whose primitive ideal space is
a circle plus
two points  (cf. \cite{hs04}) shows that no such an
algebra can be isomorphic to $C(S^2_{pq-})$. 

\section{Index computation for associated noncommutative line bundles}

\noindent
A first step in an index computation for finitely generated projective
modules is to find appropriate summable
Fredholm modules. Although the index computations for noncommutative
line bundles over the generic Podle\'s and
mirror quantum spheres proceed along the same lines (cf.~\cite{hms03}), it is 
in the structure of  Fredholm modules where  differences between these
spheres are crucial. Indeed, one can take
appropriate infinite-dimensional irreducible representations $\rho_1$
and $\rho_2$  of
$C(S^2_{pq+})$, and they yield a 1-summable Fredholm module over
$\cO(S^2_{pq+})$, the universal $*$-algebra for the defining relations
of $C(S^2_{pq+})$ \cite{mnw91,hms03}. This is in contrast with the situation 
for  mirror quantum spheres, where $(\rho_+,\rho_-)$ is not a Fredholm
module. Instead, making an appropriate identification of Hilbert 
spaces\footnote{
We are grateful to Elmar Wagner for figuring out this identification.
}, we obtain: 
\vspace*{3mm}
\ble \label{Fred-}
Let $\Tr$ be the operator trace.
The pair $(\,\rho_+\oplus\rho_-\,,\,\int_{S^1}\rho_{\lambda}d\lambda\,)$
is a 1-summable Fredholm module
over $\cO(S^2_{pq-})$, i.e.\
$
\mbox{$\fa$}\; x\in \cO(S^2_{pq-}):\;
\Tr\;\mbox{$\left|\,\rho_+(x)\oplus\rho_-(x)-\int_{S^1}\rho_{\lambda}(x)
d\lambda\,\right|$}\;
<\infty
$.
\ele
\vspace*{3mm}

On the other hand, it is straightforward to show that
$\cO(S^2_{pq-})\inc \cO(S^3_{pq\theta})$
is a 
principal extension in the sense of \cite[Definition 2.1]{bh04}.
The  $\cO(S^2_{pq-})$-modules 
associated to this extension via 1-dimensional representations
(noncommutative line bundles)
can be defined as
\beq\label{line}
\cO(S^3_{pq\theta})_{\mu}:=
\{\,x\in\cO(S^3_{pq\theta})~|~\ha_g(x)=g^{-\mu}x,\,\fa\, g\in U(1)\,\},
\quad \mu\in\Z.
\eeq
Applying \cite[Theorem 3.1]{bh04}, one can prove 
that the modules $\cO(S^3_{pq\theta})_{\mu}$ are finitely
generated projective and    can be explicitly given by the idempotents
$E_\mu=R_\mu^TL_\mu$. Here, respectively  for $\mu<0$ and  $\mu>0$,
\begin{align}
&R_\mu=(b^{|\mu|},ab^{|\mu|-1},\ldots, a^{|\mu|}),&& ~L_\mu=
\left(p^{|\mu|}A^{|\mu|}{b^*}^{|\mu|},
{\textstyle\binom{|\mu|}{1}}_p \,
p^{|\mu|-1}
A^{|\mu|-1}{b^*}^{|\mu|-1}a^*,\ldots,{a^*}^{|\mu|}\right),\\
&R_\mu=({b^*}^{\mu},a^*{b^*}^{\mu-1},\ldots, {a^*}^{\mu}),&&
~ L_\mu=
\left(b^\mu,
{\textstyle\binom{\mu}{1}}_q \,
b^{\mu-1} B a,\ldots,B^\mu a^\mu\right),
\end{align}
and  $A:=1-aa^*$, $B:=1-bb^*$. (Note that the
entries of $E_\mu$ are indeed in $\cO(S^2_{pq-})$.)
Replacing $\cO(S^3_{pq\theta})$ by $C(S^3_{pq\theta})$ in (\ref{line}),
one obtains a definition of the $C(S^2_{pq-})$-modules 
$C(S^3_{pq\theta})_\mu$ (continuous sections of noncommutative
line bundles). One can argue that these are
isomorphic to the
finitely generated projective modules $C(S^2_{pq-})^{|\mu|+1}E_\mu$
(row times matrix), with
the $E_\mu$'s given above.
\vspace*{3mm}

Now, any summable Fredholm module yields a cyclic cocycle that can be
paired with a $K$-group \cite{c-a85}. In particular, Lemma~\ref{Fred-}
gives us a trace on  $\cO(S^2_{pq-})$ defined by
$\tr(x)=\Tr((\rho_+\oplus\rho_-)(x)-\int_{S^1}\rho_{\lambda}(x)d\lambda)$.
The pairing of the cyclic 0-cocycle $\tr$ and the $K_0$-class
$[E_\mu]$, $\mu\in\Z$, yields complicated-looking rational
functions that have already been encountered in \cite{h-pm00}. Therein,
with the help of the noncommutative index formula \cite{c-a85},
they were proved to be the desired integers.
This brings us to the main
technical result:
\vspace*{3mm}
\bth\label{index}
Let $\langle~,~\rangle$ be the pairing of cyclic cohomology
and $K$-theory. Then 
$
\langle \tr,[E_\mu]\rangle = \mu$, $\fa\mu\in\Z
$.
\ethe
\vspace*{3mm}

\vspace*{3mm}\noindent
As in \cite[Corollary~2.4]{h-pm00} and 
\cite[Corollary~3.4]{hms03},
we can use  Theorem~\ref{index} to estimate the positive cone
$K_0^+(\cO(S^2_{pq-}))$.

\section{Quantum lens spaces as non-crossed-product examples of 
extensions of $C^*$-algebras}
\setcounter{equation}{0}

\noindent
Since $Z_n:=\Z/n\Z$ is a subgroup of $U(1)$, we can define new
fixed-point subalgebras
$C(L^n_{pq\theta}):=C(S^3_{pq\theta})^{Z_n}$. They play the role of
the algebras of continuous functions on quantum lens spaces
(cf.\ \cite{mt92,hs03}).
 Much as in the classical case, there is still a natural
$U(1)$-action on $C(L^n_{pq\theta})$. It can be checked that, for any
$n\in\N\setminus\{0\}$, this action is principal (free) in the sense of \cite{e-da00},
and that the fixed-point subalgebras 
$C(L^n_{pq\theta})^{U(1)}$ always coincide with $C(S^2_{pq-})$ (cf.~\cite{s-w03}).
Furthermore, we can replace $\cO(S^3_{pq\theta})$ by 
$C(L^n_{pq\theta})$ in (\ref{line}) to obtain the $C(S^2_{pq-})$-modules
$C(L^n_{pq\theta})_\mu$ and verify that they are isomorphic
with the previously defined modules:
\beq\label{nmu}
\fa\,n\in\N\setminus\{0\}:\;C(L^n_{pq\theta})_\mu\cong C(S^3_{pq\theta})_{n\mu}
\quad\mbox{as left $C(S^2_{pq-})$-modules}.
\eeq
Thus it follows from Theorem~\ref{index} that the modules 
$C(L^n_{pq\theta})_\mu$ are {\em not} free
unless $\mu=0$.
On the other hand,  using Fourier analysis with coefficients
in a  $C^*$-algebra $A$ \cite[p.222-3]{d-kr96} 
and plugging $A\!\rtimes\Z$ into 
(\ref{line}) in place of $\cO(S^3_{pq\theta})$, one can prove:
\vspace*{3mm}\ble\label{free}
Let $A$ be a unital $C^*$-algebra and $A\!\rtimes\Z$ an arbitrary 
crossed product of $A$ with $\Z$. Then, for the natural action of
$U(1)$ on $A\!\rtimes\Z$ and any $\mu\in\Z$, we have 
$(A\!\rtimes\Z)_\mu
\cong A$ as left $A$-modules.
\ele
\vspace*{3mm}
This
yields the desired contradiction, proving our main result:
\vspace*{3mm}\bth
Let $n\in \N$, $n>0$, $p,q,\theta\in [0,1)$. 
The $C^*$-algebra $C(L^n_{pq\theta})$ is a principal $U(1)$-extension 
of $C(S^2_{pq-})$, but
$C(L^n_{pq\theta})\not\cong C(S^2_{pq-})\!\rtimes\!\Z$ as $C^*$-algebras
for {\em  any} action of $\Z$ on  $C(S^2_{pq-})$.
\ethe

\vspace*{3mm}
\noindent{\bf Acknowledgements.}
\footnotesize\baselineskip9pt
The authors are very
grateful to Elmar Wagner for his help with Lemma~\ref{Fred-}, 
Sylvain Gol\'enia for the French translation,
and to Joseph C. V\'arilly for his helpful comments on the manuscript.
A preliminary version of the
full account of
this work can be found at
http://www.fuw.edu.pl/$\!\widetilde{\phantom{m}}\!$pmh.
This work was partially supported by the
European Commission
  grants MERG-CT-2004-513604 (PMH),
 MKTD-CT-2004-509794 (RM, WS),
  the KBN  grants 2 P03A 013 24 (PMH, RM, WS),
115/E-343/SPB/6.PR UE/DIE 50/2005-2008 (PMH),
and the Mid-Career Academic Grant from the Faculty of
Science and IT at the the University of Newcastle (WS).


\vspace*{-0mm}\begin{thebibliography}{00}\footnotesize\baselineskip9pt


\vspace*{-1mm}\bibitem{bhms06} P.F.~Baum, P.M.~Hajac, R.~Matthes,
W.~Szyma\'nski,
The K-theory of Heegaard-type quantum 3-spheres,
K-Theory.



\vspace*{-1mm}\bibitem{bh04}
T.~Brzezi\'nski, P.M.~Hajac,
The Chern-Galois character,
C. R. Acad. Sci. Paris, Ser. I 338 (2004) 113-116.



\vspace*{-1mm}\bibitem{c-a85} A.\  Connes,
Non-commutative differential geometry,
Inst.\ Hautes \'Etudes Sci.\ Publ.\ Math.\
        62 (1985) 257--360.

\vspace*{-1mm}\bibitem{d-kr96} K.R.~Davidson, 
$C\sp *$-algebras by example. 
Fields Institute Monographs, 6. AMS, Providence, RI, 1996.

\vspace*{-1mm}\bibitem{d-l03} L.~D\c{a}browski,
The garden of quantum spheres,  Noncommutative geometry and quantum
groups (Warsaw, 2001),  37--48, Banach Center Publ. 61, Polish Acad. Sci., Warsaw, 2003.


\vspace*{-1mm} \bibitem{dlps05} 
L.~D\c{a}browski, G.~Landi, M.~Paschke, A.~Sitarz, 
 The spectral geometry of the equatorial Podle\'s sphere.  
C. R. Math. Acad. Sci. Paris  340  (2005),   819--822.

\vspace*{-1mm}\bibitem{e-da00} D.A.~Ellwood,
A new characterisation of principal actions,
J. Funct. Anal.  173  (2000) 49--60.

\vspace*{-1mm}\bibitem{h-pm00} P.M.~Hajac,
Bundles over quantum sphere and noncommutative index theorem,
K-Theory 21 (2000) 141--150. 

\vspace*{-1mm}\bibitem{hm99} P.M.~Hajac, S.~Majid,
Projective module description of the $q$-monopole,  
Comm. Math. Phys.  206  (1999),   247--264.

\vspace*{-1mm}\bibitem{hms03} P.M.~Hajac, R.~Matthes, W.~Szyma\'nski,
Chern numbers for two families of noncommutative Hopf fibrations,
C.\ R.\ Acad.\ Sci.\  Paris Ser.~I
336 (2003) 925--930.

\vspace*{-1mm}\bibitem{hms06} P.M.~Hajac, R.~Matthes, W.~Szyma\'nski,
A locally trivial quantum Hopf fibration, Algebras Rep.\ Theory.

\vspace*{-1mm}\bibitem{hs02} J.H.~Hong, W.~Szyma\'nski,
Quantum spheres and projective spaces as graph algebras,
Commun. Math. Phys. 232 (2002) 157--188.

\vspace*{-1mm}\bibitem{hs03} J.H.~Hong, W.~Szyma\'nski,
Quantum lens spaces and graph algebras,
Pacific J. Math.  211  (2003) 249--263.

\vspace*{-1mm}\bibitem{hs04} J.H.~Hong, W.~Szyma\'nski,
The primitive ideal space of the $C^*$-algebras of infinite graphs,
J. Math. Soc. Japan 56 (2004) 45--64.

\vspace*{-1mm}\bibitem{kl93} S.~Klimek, A.~Lesniewski,
A two-parameter quantum deformation of the unit disc,
J. Funct. Anal. 115 (1993) 1--23.


\vspace*{-1mm}\bibitem{mnw91}  T.~Masuda, Y.~Nakagami, J.~Watanabe,
Noncommutative differential geometry on the quantum two sphere of Podle\'s. I: 
An algebraic viewpoint, $K$-Theory 5 (1991) 151--175.

\vspace*{-1mm}\bibitem{mt92} K.~Matsumoto, J.~Tomiyama,
Noncommutative lens spaces,
J. Math. Soc. Japan  44  (1992)
13--41.

\vspace*{-1mm}\bibitem{p-p87} P. Podle\'s,
Quantum spheres, Lett. Math. Phys. 14 (1987) 193--202.


\vspace*{-1mm}\bibitem{s-ajl91} A.J.-L.~Sheu,
Quantization of the Poisson SU(2) and its Poisson homogeneous space---the
2-sphere, Commun. Math. Phys. 135 (1991) 217--232.

\vspace*{-1mm}\bibitem{s-w03} W.~Szyma\'nski,
Quantum lens spaces and principal actions on graph $C\sp *$-algebras,
 Noncommutative geometry and quantum groups (Warsaw, 2001),  299--304,
 Banach Center Publ. 61, Polish Acad. Sci., Warsaw, 2003.

\end{thebibliography}
\end{document}